\documentclass[11pt,leqno,oneside]{article}
\usepackage{amsmath}
\usepackage{amsthm}
\usepackage{amstext}
\usepackage{amsopn}
\usepackage{amsbsy}

\newtheorem*{main-theorem}{Main Theorem}
\newtheorem{theorem}{Theorem}

\newtheorem{lemma}{Lemma}

\begin{document}

\title{A Lower Bound of The First Eigenvalue of a Closed Manifold with Positive
Ricci Curvature
}

\author{Jun LING}
\date{}

\maketitle

\begin{abstract}
We give an estimate on the lower bound of the first non-zero
eigenvalue of the Laplacian for a closed Riemannian manifold with
positive Ricci curvature in terms of the in-diameter and the lower
bound of the Ricci curvature.
\end{abstract}

\section{Introduction}\label{sec-intro}
If $M$ is an n-dimensional closed Riemannian manifold whose Ricci
curvature has a positive lower bound $(n-1)K$ for some constant
$K>0$, A.~Lichnerowicz \cite{lich} gave the following lower bound
of the first non-zero eigenvalue $\lambda$ of the Laplacian on $M$
\begin{equation}\label{nK-bound}
\lambda\geq nK.
\end{equation}
This estimate gives no information when the above constant $K$
vanishes. In such case, Li-Yau \cite{liy1} and Zhong-Yang
\cite{zy} provided another lower bound
\[
\lambda\geq \frac{\pi^2}{d^2}.
\]
It is an interesting problem to find a unified lower bound of the
first non-zero eigenvalue $\lambda$ in terms of the lower bound
$(n-1)K$ of the Ricci curvature and the diameter $d$, in-diameter
$\tilde{d}$ and other geometric quantities, which do not vanish as
$K$ vanishes, of the manifold with positive Ricci curvature.
D.~Yang \cite{yang} showed a lower bound $(1/4)(n-1)K+\pi^2/d^2$.
In this paper we give a new estimate on the lower bound of the
first non-zero eigenvalue of a closed Riemannian manifold with
positive lower bound of Ricci curvature in terms of the
in-diameter and the lower bound of Ricci curvature. Instead of
using the Zhong-Yang's canonical function or the "midrange" of the
normalized eigenfunction of the first eigenvalue in the proof,  we
use a function $\xi$ that the author constructed in \cite{ling-d}
for the construction of the suitable test function and use the
structure of the nodal domains of the eigenfunction. That provides
a new way to sharpen the bound. We have the following result.

\begin{theorem} \label{main-thm}
If $M$ is an $n$-dimensional closed Riemannian manifold. Suppose
that Ricci curvature \textup{Ric}$(M)$ of $M$ is bounded below by
$(n-1)K$ for some positive constant $K$, i. e. \textup{Ric}$(M)$
satisfies (\ref{ricci-bound})
\begin{equation}\label{ricci-bound}
\textup{Ric}(M)\geq (n-1)K,
\end{equation}
then the first non-zero eigenvalue $\lambda$ of the Laplacian of
$M$ has the following lower bound
\begin{equation}                    \label{main-bound-in}
\lambda \geq\frac12(n-1)K +  \frac{\pi^2}{\tilde{d}^2},
\end{equation}
where $\tilde{d}$ is the diameter of the largest interior ball in
the nodal domains of the first eigenfunction.
\end{theorem}

We derive some preliminary estimates and conditions for test
functions in the next section and construct the needed test
function and prove the main result in the last section.

\section{Preliminary Estimates}\label{sec-pre-es}
The classic  Lichnerowicz Theorem \cite{lich} states that if $M$
is an $n$-dimensional compact manifold without boundary whose
Ricci curvature satisfies (\ref{ricci-bound}) then the first
positive closed eigenvalue has a lower bound in (\ref{nK-bound}).
For the completeness and consistency, we use gradient estimate in
\cite{li}-\cite{liy1} and \cite{sy} to derive the Lichnerowicz
estimate.
\begin{lemma}[Lichnerowicz]               \label{nK-lemma}
Under the conditions in Theorem \ref{main-thm}, the estimate
(\ref{nK-bound}) holds.
\end{lemma}
\begin{proof}\quad
Let $v$ be a normalized eigenfunction of the first closed
eigenvalue such that
\begin{equation}                                \label{v-c-con}
\sup_{M}v=1, \quad \inf_{M}v=-k
\end{equation}
with $0<k\leq 1$. The function $v$ satisfies the following
equation
\begin{equation}                                 \label{v-c-eq}
\Delta v=-\lambda v\quad \textrm{in }M,
\end{equation}
where $\Delta$ is the Laplacian of $M$.

Take an orthonormal frame $\{e_1, \dots, e_n\}$ of $M$ about
$x_0\in M$. At $x_0$ we have
\[
\nabla_{e_j}(|\nabla v|^2)(x_0)= \sum_{i=1}^n2v_iv_{ij}
\]
and
\begin{eqnarray}
\Delta (|\nabla v|^2)(x_0)&=&2\sum_{i,j=1}^nv_{ij}v_{ij}+2\sum_{i,j=1}^nv_{i}v_{ijj}\nonumber\\
      &=&2\sum_{i,j=1}^nv_{ij}v_{ij}+2\sum_{i,j=1}^nv_{i}v_{jji}+2\sum_{i,j=1}^n\textrm{R}_{ij}v_iv_j\nonumber\\
      &=&2\sum_{i,j=1}^nv_{ij}v_{ij}+2\nabla v \nabla (\Delta v) + 2\textrm{Ric}(\nabla v, \nabla v)\nonumber\\
      &\geq&2\sum_{i=1}^nv_{ii}^2+ 2\nabla v \nabla (\Delta v) + 2(n-1)K|\nabla
      v|^2\nonumber\\
      &\geq&\frac{2}{n}(\Delta v)^2-2\lambda|\nabla v |^2 + 2(n-1)K|\nabla v|^2.\nonumber
\end{eqnarray}
Thus at all point $x\in M$,
\begin{equation}                        \label{basic2.002-c}
\frac12\Delta (|\nabla v|^2)\geq \frac1n \lambda^2 v^2 +
[(n-1)K-\lambda]|\nabla v|^2.
\end{equation}
On the other hand, after multiplying (\ref{v-c-eq}) by $v$ and
integrating both sides over $M$, we have
\[
\int_{M}\lambda v^2\,dx=-\int_Mv\Delta v\, dx =\int_{M}|\nabla
v|^2\, dx.
\]
Integrating (\ref{basic2.002-c}) over $M$ and using the above
equality, we get
\begin{equation}                                        \label{basic2.003-c}
0\geq\int_M(nK-\lambda)\frac{n-1}{n}\lambda v^2\,dx.
\end{equation}
Therefore (\ref{nK-bound}) holds.
\end{proof}

\begin{lemma}       \label{pre-es-c}
Let $v$ be, as the above, the normalized eigenfunction for the
first non-zero eigenvalue $\lambda$. Then $v$ satisfies the
following
\begin{equation}                        \label{basic3-c}
\frac{\left |\nabla v\right |^2}{b^2-v^2} \leq\lambda,
\end{equation}
where $b>1$ is an arbitrary constant.
\end{lemma}
\begin{proof}\quad
Consider the function
\begin{equation}                \label{p-of-x-def}
P(x)=|\nabla v|^2+Av^2,
\end{equation}
where $A=\lambda (1 +\epsilon)$ for small $\epsilon>0$. Function
$P$ must achieve its maximum at some point $x_0\in M$. We claim
that $\nabla v (x_0)=0$.

If on the contrary, $\nabla v (x_0)\not=0$, then we can rotate the
local orthonormal about $x_0$ such that
\[
|v_1(x_0)|=|\nabla v(x_0)|\not=0\qquad \textrm{and}\qquad
v_i(x_0)=0,\quad i\geq2.
\]
Since $P$ achieves its maximum at $x_0$, we have,
\[
\nabla P(x_0)=0 \qquad \textrm{and}\qquad \Delta P(x_0)\leq 0.
\]
That is, at $x_0$ we have
\[
0=\frac12\nabla_{i} P=\sum_{j=1}^{n}v_jv_{ji}+Avv_i,
\]
\begin{equation}                \label{d1}
v_{11}=-Av \qquad \textrm{and}\qquad v_{1i}=0\quad i\geq 2,
\end{equation}
and
\begin{eqnarray}
&0 &\geq \frac12\Delta P(x_0) =\sum_{i, j=1}^{n}\left(v_{ji}v_{ji}+v_{j}v_{jii}+Av_{i}v_{i}+Avv_{ii}\right)\nonumber\\
&{} &=\sum_{i, j=1}^{n}\left(v_{ji}^2+v_j(v_{ii})_j +\textrm{R}_{ji}v_{j}v_{i}+Av_{ii}^2 +Av v_{ii}\right)\nonumber\\
&{} &=\sum_{i, j=1}^{n}v_{ji}^2+\nabla v\nabla(\Delta v) +\textrm{Ric}(\nabla v,\nabla v)+A|\nabla v|^2 +Av\Delta v\nonumber\\
&{} &\geq v_{11}^2+\nabla v\nabla(\Delta v) + (n-1)K|\nabla v|^2+A|\nabla v|^2 +Av\Delta v\nonumber\\
&{} &=(-Av)^2-\lambda |\nabla v|^2+ (n-1)K|\nabla v|^2+A|\nabla v|^2 -\lambda Av^2\nonumber\\
&{} &=(A-\lambda + (n-1)K)|\nabla v|^2+Av^2(A-\lambda)\nonumber
\end{eqnarray}
where we have used (\ref{d1}) and (\ref{ricci-bound}). Therefore
at $x_0$,
\begin{equation}
0\geq (A-\lambda)|\nabla v|^2+A(A-\lambda)v^2, \label{d2}
\end{equation}
that is,
\[
|\nabla v(x_0)|^2+ \lambda(1+\epsilon)v(x_0)^2\leq 0.
\]
Thus $\nabla v(x_0) = 0$. This contradicts $\nabla v(x_0)\not= 0$.
So the above claim is right.

Therefore  we have $\nabla v(x_0)=0$,
\[
P(x_0)=|\nabla v(x_0)|^2 +Av(x_0)^2=Av(x_0)^2\leq A,
\]
and at all $x\in M$
\[
|\nabla v(x)|^2+Av(x)^2=P(x)\leq P(x_0)\leq A.
\]
Letting $\epsilon \rightarrow 0$ in the above inequality, the
estimate (\ref{basic3-c}) follows.
\end{proof}

We want to improve the upper bound in (\ref{basic3-c}) further and
proceed in the following way.

Define a function $Z$ on $[-\sin^{-1}(k/b),\sin^{-1}(1/b)]$ by
\[
Z(t)=\max_{x\in M,t=\sin^{-1} \left(v\left(x\right)/b\right)}
\frac{\left |\nabla v\right |^2}{b^2-v^2}/\lambda.
\]
From (\ref{basic3-c}) we have
\begin{equation}                                \label{basic5-c}
Z(t)\leq 1\qquad \textrm{on } [-\sin^{-1}(k/b),\sin^{-1}(1/b)]
\end{equation}
For convenience, in this paper we let
\begin{equation}\label{alpha-delta}
\alpha =
\frac12(n-1)K\qquad\textrm{and}\qquad\delta=\alpha/\lambda.
\end{equation}
By (\ref{nK-bound}) we have
\begin{equation}\label{delta-bound}
\delta\leq \frac{n-1}{2n}.
\end{equation}
We have the following conditions on the test function $Z$.

\begin{theorem}                                                     \label{thm-barrier-c}
If the function $z:\ [-\sin^{-1}(k/b),\sin^{-1}(1/b)]\mapsto
\mathbf{R}^1$ satisfies the following
\begin{enumerate}
 \item $z(t)\geq Z(t) \qquad t\in [-\sin^{-1}(k/b),\sin^{-1}(1/b)]$,
 \item there exists some $x_0\in M$
       such that at point $t_0=\sin^{-1} (v(x_0)/b)$ \linebreak
       $z(t_0)=Z(t_0)$,
 \item $z(t_0)>0$,\qquad\textrm{and}
 \item $z'(t_0)\sin t_0\geq 0$,
\end{enumerate}
then we have the following
\begin{equation}                            \label{barrier-eq-c}
0\leq\frac12z''(t_0)\cos^2t_0 -z'(t_0)\cos t_0\sin t_0 -z(t_0)+
1-2\delta \cos^2t_0.
\end{equation}
\end{theorem}

\begin{proof}\quad
Define
\[ J(x)=\left\{ \frac{\left |\nabla v\right |^2}{b^2-v^2}
-\lambda z \right\}\cos^2t,
\]
where $t=\sin^{-1}(v(x)/b)$. Then
\[ J(x)\leq 0\quad\textrm{for } x\in M
\qquad \textrm{and} \qquad J(x_0)=0.
\]
If $\nabla v(x_0)=0$, then
\[ 0=J(x_0)=-\lambda z\cos^2 t.
\]
This contradicts Condition 3 in the theorem. Therefore
\[ \nabla v(x_0)\not=0.
\]

The Maximum Principle implies that
\begin{equation}                                                \label{es1}
\nabla J(x_0)=0\qquad \textrm{and}\qquad \Delta J(x_0)\leq 0.
\end{equation}
$J(x)$ can be rewritten as
\[  J(x)=\frac{1}{b^2}|\nabla v|^2-\lambda z\cos^2t.
\]
Thus (\ref{es1}) is equivalent to
\begin{equation}                                              \label{es2}
\frac{2}{b^2}\sum_{i}v_iv_{ij}\Big|_{x_0}=\lambda\cos t[z' \cos t
-2z\sin t]t_j\Big|_{x_0}
\end{equation}
and
\begin{eqnarray}                                              \label{es3}
0&\geq&\frac{2}{b^2}\sum_{i,j}v_{ij}^2+\frac{2}{b^2}\sum_{i,j}v_iv_{ijj}
 -\lambda (z''|\nabla t|^2+z'\Delta t)\cos^2t \\
 & &+4\lambda z'\cos t\sin t |\nabla t|^2 -
\lambda z\Delta\cos^2t\Big|_{x_0}.\nonumber
\end{eqnarray}
Rotate the frame so that $v_1(x_0)\not=0$ and $v_i(x_0)=0$ for
$i\geq 2$. Then (\ref{es2}) implies
\begin{equation}                                             \label{es4}
v_{11}\Big|_{x_0}=\frac{\lambda b}{2}(z'\cos t-2z\sin t)
\Big|_{x_0}\quad\text{and}\quad v_{1i} \Big|_{x_0}=0\ \text{for }
i\geq2.
\end{equation}
Now we have
\begin{eqnarray}
|\nabla v|^2
\Big|_{x_0}&=&\lambda b^2z\cos^2t\Big|_{x_0},\nonumber\\
 |\nabla t|^2
\Big|_{x_0}&=&\frac{|\nabla v|^2}{b^2-v^2}=\lambda z
\Big|_{x_0},\nonumber\\
\frac{\Delta v}{b}\Big|_{x_0} &=&\Delta \sin t =\cos t\Delta
t-\sin t |\nabla t|^2
\Big|_{x_0},\nonumber\\
\Delta t\Big|_{x_0}&=&\frac{1}{\cos t}(\sin t|\nabla
t|^2+\frac{\Delta v}{b})
\nonumber\\
 &=&\frac{1}{\cos t} [ \lambda z\sin t-\frac{\lambda}{b}v] \Big|_{x_0}, \quad\textrm{and}
\nonumber\\
\Delta\cos^2t\Big|_{x_0}&=&\Delta \left(1-\frac{v^2}{b^2}\right)
 =-\frac{2}{b^2}|\nabla v|^2-\frac{2}{b^2}v\Delta v
\nonumber\\
  &=&-2\lambda z\cos^2t+\frac{2}{b^2}\lambda v^2\Big|_{x_0}. \nonumber
\end{eqnarray}
Therefore,
\begin{eqnarray}
& {}&
\frac{2}{b^2}\sum_{i,j}v_{ij}^2\Big|_{x_0}\geq\frac{2}{b^2}v_{11}^2
\nonumber\\
& {}& =\frac{\lambda ^2}{2}(z')^2\cos^2t-2\lambda ^2zz'\cos t\sin
t
      +2\lambda ^2z^2\sin^2 t\Big|_{x_0}\nonumber,
\end{eqnarray}
\begin{eqnarray}
\frac{2}{b^2}\sum_{i,j}v_iv_{ijj}\Big|_{x_0}
&=&\frac{2}{b^2}\left(\nabla v\,\nabla
       (\Delta v)+\textrm{Ric}(\nabla v,\nabla v)\right)\nonumber\\
&\geq& \frac{2}{b^2}(\nabla v\,\nabla (\Delta v)+(n-1)K|\nabla v|^2)\nonumber\\
 &=&-2\lambda^2z\cos^2t+4\alpha \lambda z\cos^2t\Big|_{x_0},\nonumber
\end{eqnarray}
\begin{eqnarray}
&{}&  -\lambda (z''|\nabla t|^2+ z'\Delta t)\cos^2t\Big|_{x_0}\nonumber\\
&{}&=-\lambda^2 zz''\cos^2t-
\lambda^2zz'\cos t\sin t\nonumber\\
&{}&{ }+\frac{1}{b}\lambda^2z'v\cos t \Big|_{x_0},\nonumber
\end{eqnarray}
and
\begin{eqnarray}
&{}&4\lambda z'\cos t\sin t|\nabla t|^2-\lambda z\Delta
\cos^2t\Big|_{x_0}
\nonumber\\
&{}&=4\lambda^2zz'\cos t\sin
t+2\lambda^2z^2\cos^2t-\frac{2}{b}\lambda^2zv\sin t
\Big|_{x_0}.\nonumber
\end{eqnarray}
Putting these results into (\ref{es3}) we get
\begin{eqnarray}                                                   \label{es5}
0&\geq&-\lambda^2zz''\cos^2t+ \frac{\lambda^2}{2}(z')^2\cos^2t
+\lambda^2z'\cos t\left(z\sin t +\sin t\right)
  \nonumber\\
 & & {}+2\lambda^2z^2-2\lambda^2z +4\alpha \lambda z\cos^2t
 \Big|_{x_0},
\end{eqnarray}
where we used (\ref{es4}). Now
\begin{equation}                                                    \label{es6}
z(t_0)>0,
\end{equation}
by Condition 3 in the theorem. Dividing two sides of (\ref{es5})
by $2\lambda^2z\Big|_{x_0}$, we have
\begin{eqnarray}
0&\geq&-\frac12z''(t_0)\cos^2t_0 +\frac12z'(t_0)\cos t_0\left(\sin t_0+\frac{\sin t_0}{z(t_0)}\right) +z(t_0) \nonumber\\
 & & {}  -1 +2\delta \cos^2t_0+\frac{1}{4z(t_0)}(z'(t_0))^2\cos^2t_0.\nonumber
\end{eqnarray}
Therefore,
\begin{eqnarray}
0&\geq&-\frac12z''(t_0)\cos^2t_0 + z'(t_0)\cos t_0\sin t_0+z(t_0) -1 +2\delta \cos^2t_0\nonumber\\
 & & {}+\frac{1}{4z(t_0)}(z'(t_0))^2\cos^2t_0+\frac12z'(t_0)\sin t_0\cos
 t_0[\frac{1}{z(t_0)}-1].                               \label{es6.1}
\end{eqnarray}
Conditions 1, 2 and 4 in the theorem imply that
$0<z(t_0)=Z(t_0)\leq 1$ and $z'(t_0)\sin t_0\geq 0$. Thus the last
two terms in (\ref{es6.1}) are nonnegative and
(\ref{barrier-eq-c}) follows.
\end{proof}

\section{Proof of the Main Result}\label{sec-proof}

\begin{proof}[Proof of Theorem \ref{main-thm}]\quad
Let
\begin{equation}                                \label{z-def}
z(t)=1+\delta\xi(t),
\end{equation}
where $\xi$ is the functions defined by (\ref{xi-def}) in Lemma
\ref{xi-lemma}. We claim that
\begin{equation}                        \label{4.1-c}
Z(t)\leq z(t)\qquad \textrm{on } [-\sin^{-1}(k/b),\sin^{-1}(1/b)].
\end{equation}

Lemma \ref{xi-lemma} implies that for $t\in
[-\sin^{-1}(k/b),\sin^{-1}(1/b)] $ we have the following
\begin{eqnarray}
& &{}\frac{1}{2}z''\cos ^2t-z'\cos t\sin t-z
    =-1+ 2\delta\cos^2t,               \label{z-eq}\\
& &{}z'(t)\sin t\geq 0, \label{z'-geq0}\\
& &{}0<1-(\frac{\pi^2}{4}-1)\frac{n-1}{2n}\leq
1-(\frac{\pi^2}{4}-1)\delta=z(0)\leq z(t),\quad\textrm{and}\label{z-min}\\
& &{}z(t) \leq z(\frac{\pi}{2})=1. \label{z-max}
\end{eqnarray}

Let $P\in\mathbf{R}^1$ and $t_0\in
[-\sin^{-1}(k/b),\sin^{-1}(1/b)]$ such that
\[ P=\max_{t\in [-\sin^{-1}(k/b),\sin^{-1}(1/b)]}\left(Z(t)-z(t)\right)=Z(t_0)-z(t_0).
\]
Thus
\[
Z(t)\leq z(t)+P\quad \textrm{for }t\in
[-\sin^{-1}(k/b),\sin^{-1}(1/b)]\quad\textrm{and}\quad
Z(t_0)=z(t_0)+P.
\]
Suppose that $P>0$. Then $z+P$ satisfies the conditions in Theorem
\ref{thm-barrier-c}. (\ref{barrier-eq-c}) implies
\begin{eqnarray}
&{}&z(t_0)+P=Z(t_0)\nonumber\\
&{}&\leq  \frac12(z+P)''(t_0)\cos^2 t_0-(z+P)'(t_0)\cos t_0
 \sin t_0+1-2\delta \cos^2 t_0\nonumber\\
&{}&=\frac12z''(t_0)\cos^2t_0-z'(t_0)\cos t_0\sin
t_0+1-2\delta \cos^2 t_0\nonumber\\
&{}&=z(t_0).\nonumber
\end{eqnarray}
This contradicts the assumption $P>0$. Thus $P\leq 0$ and
(\ref{4.1-c}) holds. That means
\begin{equation}\label{4.3-c}
\sqrt{\lambda}\geq \frac{|\nabla t|}{\sqrt{z(t)}}.
\end{equation}

Note that the eigenfunction $v$ of the first nonzero eigenvalue
has exactly two nodal domains $D^+=\{x: v(x)>0\}$ and $D^-=\{x:
v(x)<0\}$ and the nodal set $v^{-1}(0)$ is compact (see \cite{cha}
and \cite{ch}). Take $q_1$ on $M$ such that $v(q_1)=1 =\sup_M v$
and and $q_2\in v^{-1}(0)$ such that distance $d(q_1, q_2) =
\textrm{ distance }d(q_1, v^{-1}(0))$. Let $L$ be the minimum
geodesic segment between $q_1$ and $q_2$. We integrate both sides
of (\ref{4.3-c}) along $L$ and change variable and let
$b\rightarrow 1$. Let $d_+$, $d_-$ be the diameter of the largest
interior ball in $D^+$, $D^-$ respectively,
\[
d_+=2r_+ \qquad\textup{and}\qquad r_+=\max_{x\in D_+}
\textup{dist}(x, v^{-1}(0))
\]
and
\[
d_-=2r_- \qquad\textup{and}\qquad r_-=\max_{x\in D_-}
\textup{dist}(x, v^{-1}(0)).
\]
 Then $\tilde{d}=\max\{d_+, d_-\}$ and
\begin{equation}\label{4.4-c}
\sqrt{\lambda}\,\frac{d_+}{2}\geq \int_{L}\, \frac{|\nabla
t|}{\sqrt{z(t)}}\, dl=\int_0^{\frac{\pi}{2}}
\frac{1}{\sqrt{z(t)}}\,dt \geq \frac{\left(\int_0^{\pi/2}\
\,dt\right)^\frac32}{(\int_0^{\pi/2}\ z(t)\,dt)^{\frac12}} \geq
\left( \frac{(\frac{\pi}{2})^3}{\int_0^{\pi/2}\  z(t)\,dt}
\right)^{\frac12}
\end{equation}
Square the two sides. Then
\[
\lambda \geq \frac{\pi^3}{2(d_+)^2\int_0^{\pi/2} \ z(t)\,dt}.
\]
Now
\[
\int_0^{\frac{\pi}{2}}\ z(t)\,dt=\int_0^{\frac{\pi}{2}}\ [1+
\delta \xi(t)]\,dt=\frac{\pi}{2}(1-\delta),
\]
by (\ref{xi-int}) in Lemma \ref{xi-lemma}. That is,
\[
\lambda \geq \frac{\pi^2}{(1-\delta)(d_+)^2}\quad\textrm{and}\quad
\lambda \geq  \frac12(n-1)K+\frac{\pi^2}{(d_+)^2}.
\]
Noticing that $\tilde{d}\geq d_+$ and $\tilde{d}\geq d_-$, we
complete the proof.
\end{proof}

We now present a lemma that is used in the proof of Theorem
\ref{main-thm}.

\begin{lemma}                           \label{xi-lemma}
Let
\begin{equation}                        \label{xi-def}
\xi(t)=\frac{\cos^2t+2t\sin t\cos t +t^2-\frac{\pi^2}{4}}{\cos^2t}
\qquad \textrm{on}\quad [-\frac{\pi}{2},\frac{\pi}{2}\,].
\end{equation}
Then the function $\xi$  satisfies the following
\begin{eqnarray}
& &{}\frac{1}{2}\xi''\cos ^2t-\xi'\cos t\sin t-\xi
    =2\cos^2t\quad \textrm{in }(-\frac{\pi}{2},\frac{\pi}{2}\,),          \label{xi-eq}\\
& &{}\xi'\cos t -2\xi\sin t =4t\cos t                      \label{xi-eq2}\\
& &{}\int_0^{\frac{\pi}{2}}\xi(t)\, dt= -\frac{\pi}{2}          \label{xi-int}\\
& &{}1-\frac{\pi^2}{4}=\xi(0)\leq\xi(t)\leq\xi(\pm
\frac{\pi}{2})=0\quad
\textrm{on }[-\frac{\pi}{2},\frac{\pi}{2}\,],                                \nonumber\\
& &{} \xi' \textrm{ is increasing on }
[-\frac{\pi}{2},\frac{\pi}{2}\,] \textrm{ and }
\xi'(\pm \frac{\pi}{2}) =\pm \frac{2\pi}{3},                     \nonumber\\
& &{}\xi'(t)< 0 \textrm{ on }(-\frac{\pi}{2},0)\textrm{ and \ }
\xi'(t)>0 \textrm{ on }(0,\frac{\pi}{2}\,), \nonumber \\
& &{}\xi''(\pm\frac{\pi}{2})=2, \ \xi''(0)=2(3-\frac{\pi^2}{4})
\textrm{ and \ } \xi''(t)> 0 \textrm{ on }
[-\frac{\pi}{2},\frac{\pi}{2}\,], \nonumber\\
& &{}(\frac{\xi'(t)}{t})'>0 \textrm{ on } (0,\pi/2\,)\textrm{ and
\ } 2(3-\frac{\pi^2}{4})\leq \frac{\xi'(t)}{t}\leq \frac43
\textrm{ on } [-\frac{\pi}{2},\frac{\pi}{2}\,],\nonumber \\
& &{}\xi'''(\frac{\pi}{2})=\frac{8\pi}{15}, \xi'''(t)< 0 \textrm{
on }(-\frac{\pi}{2},0) \textrm{ and \ }  \xi'''(t)>0 \textrm{ on
}(0,\frac{\pi}{2}\,). \nonumber
\end{eqnarray}
\end{lemma}
\begin{proof}\quad
For convenience, let $q(t)= \xi'(t)$, i.e.,
\begin{equation}                                       \label{q-def}
q(t) = \xi'(t) = \frac{2(2t\cos t +t^2\sin t +\cos^2 t \sin t
-\frac{\pi^2}{4}\sin t)}{\cos^3 t}.
\end{equation}
Equation (\ref{xi-eq}) and the values $\xi(\pm \frac{\pi}{2})=0$,
$\xi(0)=1-\frac{\pi^2}{4}$ and $\xi'(\pm \frac{\pi}{2}) =\pm
\frac{2\pi}{3}$ can be verified directly from (\ref{xi-def}) and
(\ref{q-def}) .  The values of $\xi''$ at $0$ and $\pm
\frac{\pi}{2}$ can be computed via (\ref{xi-eq}). By
(\ref{xi-eq2}), $(\xi(t)\cos^2 t)' =4t\cos^2 t$. Therefore
\newline $\xi(t)\cos^2 t=\int_{\frac{\pi}{2}}^t \ 4s\cos^2 s\,ds$,
and
\[
\int_{-\frac{\pi}{2}}^{\frac{\pi}{2}}\
\xi(t)\,dt=2\int_0^{\frac{\pi}{2}}\
\xi(t)\,dt=-8\int_0^{\frac{\pi}{2}}\left( \frac{1}{\cos^2(t)}
\int_t^{\frac{\pi}{2}}\ s\cos^2s\,ds\right)\,dt
\]
\[
=-8\int_0^{\frac{\pi}{2}}\left(\int_0^s\
\frac{1}{\cos^2(t)}\,dt\right)\ s\cos^2s\,ds
=-8\int_0^{\frac{\pi}{2}}\ s\cos s\sin s\,ds=-\pi.
\]
It is easy to see that $q$ and $q'$ satisfy the following
equations
\begin{equation}                                         \label{q-eq}
\frac12 q''\cos t -2q'\sin t -2q\cos t = -4 \sin t,
\end{equation}
and
\begin{equation}                                       \label{q'-eq}
\frac{\cos^2 t}{2(1+\cos^2 t)}(q')''-\frac{2\cos t\sin t}{1+\cos^2
t}(q')'-2(q')=-\frac{4}{1+\cos^2 t}.
\end{equation}
The last equation implies $q'=\xi''$ cannot achieve its
non-positive local minimum at a point in $(-\frac{\pi}{2},
\frac{\pi}{2})$. On the other hand, $\xi''(\pm\frac{\pi}{2})=2$,
by equation (\ref{xi-eq}), $\xi(\pm \frac{\pi}{2})=0$ and
$\xi'(\pm \frac{\pi}{2})=\pm \frac{2\pi}{3}$. Therefore
$\xi''(t)>0$ on $[-\frac{\pi}{2},\frac{\pi}{2}]$ and $\xi'$ is
increasing. Since $\xi'(t)=0$, we have $\xi'(t)< 0$ on
$(-\frac{\pi}{2},0)$ and $\xi'(t)>0$ on $(0,\frac{\pi}{2})$.
Similarly, from the equation
\begin{eqnarray}                                           \label{q''-eq}
&\frac{\cos^2 t}{2(1+\cos^2 t)}(q'')'' -\frac{\cos t\sin t
(3+2\cos^2 t)}{(1+\cos^2 t)^2}(q'')' -\frac{2(5\cos^2 t+\cos^4 t)}{(1+\cos^2 t)^2}(q'') \nonumber\\
&=-\frac{8\cos t\sin t}{(1+\cos^2 t)^2}
\end{eqnarray}
we get the results in the last line of the lemma.

Set $h(t)=\xi''(t)t-\xi'(t)$. Then $h(0)=0$  and $h'(t)=
\xi'''(t)t>0$ in $(0,\frac{\pi}{2})$. Therefore
$(\frac{\xi'(t)}{t})'=\frac{h(t)}{t^2}>0$ in $(0,\frac{\pi}{2})$.
Note that $\frac{\xi'(-t)}{-t}= \frac{\xi'(t)}{t}$,
$\frac{\xi'(t)}{t}|_{t=0}=\xi''(0)=2(3-\frac{\pi^2}{4})$ and
$\frac{\xi'(t)}{t}|_{t=\pi/2}=\frac43$. This completes the proof
of the lemma.
\end{proof}

%
%

Department of mathematics, Utah Valley State College, Orem, Utah 84058

\textit {E-mail address}: \texttt{lingju@uvsc.edu}
\end{document}